\documentclass{article}
\usepackage{amssymb}

\usepackage{amsmath}

\newtheorem{theorem}{Theorem}

\numberwithin{theorem}{section}
\numberwithin{equation}{section}

\begin{document}

\title{\textbf{A short proof of chaos in an atmospheric system}}
\author{{ Petre Birtea$^{1}$, Mircea Puta$^{2}$, Tudor
S. Ratiu$^{1}$, R\u azvan Tudoran$^{1}$}}
\addtocounter{footnote}{1}
\footnotetext{Institut Bernoulli, 
\'Ecole Polytechnique F\'ed\'erale de Lausanne. CH--1015 Lausanne.
Switzerland.}
\addtocounter{footnote}{1}
\footnotetext{Departamentul de Matematic\u a, Universitatea de Vest, 1900
Timi\c soara, Romania. }
\date{March 8, 2002}
\maketitle

\begin{abstract}
We will prove the presence of chaotic motion in the Lorenz five-component atmospheric system model
using the Melnikov function method developed by
Holmes and Marsden for Hamiltonian systems on Lie Groups.
\medskip

\medskip

\end{abstract}

\section{Introduction.}

The first model equations for the atmosphere are the so called primitive
equations (PE). This model allows wave-like motions on different time
scales. The slow motions which have a period of order of days are called
Rossby waves and the fast motions which have a period of hours are called
gravity waves. The question of how to balance these two time scales lead
Lorenz \cite{lor} to introduce a simplified version of the (PE) model, the
so called five-component model. This is a system of five differential
equations which couples the Rossby waves and gravity waves. This system
turns out to have a Poisson formulation on $\mathbb{R}^{5}$ first discovered by Bokhove \cite{bok}. We shall find a Poisson diffeomorphism between the Poisson structure of Bokhove and the product structure on $\mathfrak{se}^{\ast}(2)\times \mathbb{R}^{2}$, where the first factor is endowed with the Lie-Poisson structure  and the second with the standard symplectic structure. Using this diffeomorphism, the five-component Lorenz model takes on a form appropriate for the application of the Melnikov function method on Lie groups developed by Holmes and Marsden \cite{H-M}.

\section{The geometry of Lorenz simplified 
model of Rossby-gravity wave
interaction.}

The model introduced by Lorenz in \cite{lor}
is described by the following
set of differential equations:
\begin{eqnarray}
&&
\begin{array}{cc}
\overset{\cdot }{x}_{1}= 
& -x_{2}x_{3}+\varepsilon x_{2}x_{5}
\end{array}
\notag \\
&&
\begin{array}{cc}
\overset{\cdot }{x}_{2}= & x_{1}x_{3}
-\varepsilon x_{1}x_{5}
\end{array}
\notag \\
&&
\begin{array}{cc}
\overset{\cdot }{x}_{3}= & -x_{1}x_{2}
\end{array}
\label{sis1} \\
&&
\begin{array}{cc}
\overset{\cdot }{x}_{4}= & -x_{5}
\end{array}
\notag \\
&&
\begin{array}{cc}
\overset{\cdot }{x}_{5}= 
& x_{4}+\varepsilon x_{1}x_{2},
\end{array}
\notag
\end{eqnarray}%
where the variables $x_{4},$ $x_{5}$\ represent the fast gravity wave
oscillations and$\ x_{1}$, $x_{2}$, $x_{3}$\ are the slow Rossby wave
oscillations, with the parameter$\ \varepsilon $ (related to the Rossby
number) coupling the two sets of variables.

In \cite{bok}, Bokhove wrote the system (\ref{sis1}) in the following
Hamiltonian form:
\begin{equation*}
\overset{\cdot }{x}= \{x,H\}_1,
\end{equation*}
where the Hamiltonian function is given by
\begin{equation*}
H(x_{1},x_{2},x_3, x_4 ,x_{5})=\frac{1}{2}
(x_{1}^{2}+2x_{2}^{2}+x_{3}^{2}+x_{4}^{2}+x_{5}^{2}).
\end{equation*}
and the Poisson bracket $\{\cdot ,\cdot \}_{1}$ is defined by
\begin{eqnarray*}
\{f,g\}_{1} &=&
x_{1} \left[ \left( 
\frac{\partial f}{\partial x_{2}}
\frac{\partial g}{\partial x_{3}}
-\frac{\partial f}{\partial x_{3}}
\frac{\partial g}{\partial x_{2}}\right)
+\varepsilon \left( \frac{\partial f}{\partial x_{5}}
\frac{\partial g}{\partial x_{2}}
-\frac{\partial f}{\partial x_{2}}
\frac{\partial g}{\partial x_{5}}
\right)\right] \\
&&\quad + x_{2} 
\left[ \left(
\frac{\partial f}{\partial x_{3}}
\frac{\partial g}{\partial x_{1}}-
\frac{\partial f}{\partial x_{1}}
\frac{\partial g}{\partial x_{3}}
\right)+\varepsilon 
\left(\frac{\partial f}{\partial x_{1}}
\frac{\partial g}{\partial x_{5}}
-\frac{\partial f}{\partial x_{5}}
\frac{\partial g}{\partial x_{1}}
\right)\right]
\\
&& \quad +
\frac{\partial f}{\partial x_{5}}
\frac{\partial g}{\partial x_{4}}-
\frac{\partial f}{\partial x_{4}}
\frac{\partial g}{\partial x_{5}}.
\end{eqnarray*}

On the space $\mathfrak{se}^{\ast }(2)\times \mathbb{R}^{2}$ consider the Poisson
bracket given by the product bracket of the Lie-Poisson bracket on $\mathfrak{se}^{\ast }(2)$ and the Poisson bracket on $\mathbb{R}^{2}$ induced by the
standard symplectic form on $\mathbb{R}^{2}$. Denoting by $(\mu _{1},\mu _{2},\mu _{3})$ the variables on $\mathfrak{se}^{\ast }(2)$
and by $(u_{1},u_{2})$ the variables on $\mathbb{R}^{2}$, the product Poisson bracket is given by 
\begin{align*}
\{f,g\}_{2}
=&\mu _{1}\left(\frac{\partial f}{\partial \mu _{2}}
\frac{\partial g}{\partial \mu _{3}}
-\frac{\partial f}{\partial \mu _{3}}
\frac{\partial g}{\partial \mu _{2}}\right)
+\mu _{2}\left(\frac{\partial f}{\partial \mu _{3}}
\frac{\partial g}{\partial \mu _{1}}
-\frac{\partial f}{\partial \mu _{1}}
\frac{\partial g}{\partial \mu _{3}}\right) \\
& \quad+\frac{\partial f}{\partial u_{2}}
\frac{\partial g}{\partial u_{1}}
-\frac{\partial f}{\partial u_{1}}
\frac{\partial g}{\partial u_{2}}\,.
\end{align*}
The Casimir function for this
Poisson bracket is given by
\begin{equation*}
C(\mu _{1},\mu _{2},\mu _{3})=\mu _{1}^{2}+\mu _{2}^{2}.
\end{equation*}

It is easy to verify that the linear transformation 
$\Phi :\mathbb{R}^{5}\rightarrow \mathfrak{se}^{\ast }(2)\times \mathbb{R}^{2}$ given by
\begin{equation*}
\Phi (x_{1},x_2, x_3, x_4, x_{5}
=(x_{1},x_{2},x_{3},x_{4},\varepsilon x_{3}+x_{5}).
\end{equation*}
is a Poisson diffeomorphism between 
$(\mathbb{R}^{5},\{\cdot ,\cdot \}_{1})$
and $(\mathfrak{se}^{\ast }(2)\times \mathbb{R}^{2},\{\cdot ,\cdot \}_{2})$.

In the new variables the system (\ref{sis1}) becomes
\begin{eqnarray}
&&
\begin{array}{cc}
\overset{\cdot }{\mu }_{1}
= & -\mu _{2}\mu _{3}+\varepsilon \mu
_{2}u_{2}-\varepsilon ^{2}\mu _{2}\mu _{3}
\end{array}
\notag \\
&&
\begin{array}{cc}
\overset{\cdot }{\mu }_{2}
= & \mu _{1}\mu _{3}-\varepsilon \mu
_{1}u_{2}+\varepsilon ^{2}\mu _{1}\mu _{3}
\end{array}
\notag \\
&&
\begin{array}{cc}
\overset{\cdot }{\mu }_{3}= & -\mu _{1}\mu _{2}
\end{array}
\label{sis2} \\
&&
\begin{array}{cc}
\overset{\cdot }{u}_{1}= & -u_{2}+\varepsilon \mu _{3}
\end{array}
\notag \\
&&
\begin{array}{cc}
\overset{\cdot }{u}_{2}= & u_{1}, 
\end{array}
 \notag
\end{eqnarray}
which is an Hamiltonian system with respect to the Poisson bracket $\{\cdot ,\cdot \}_{2}$ and the Hamiltonian function is given by
\begin{equation}
H^{\varepsilon }(\mu _{1},\mu _{2},\mu _{3},u_{1},u_{2})=\frac{1}{2}(\mu_{1}^{2}+2\mu _{2}^{2}
+\mu _{3}^{2}+u_{1}^{2}+u_{2}^{2}
-2\varepsilon \mu_{3}u_{2}+\varepsilon ^{2}\mu _{3}^{2}).  \label{hap}
\end{equation}

\section{Chaos by the Melnikov method.}

We will prove the occurrence of chaotic motion in the system (\ref{sis2}) by
showing the existence of transverse heteroclinic orbits. Melnikov \cite{mel}
gave an effective method to prove the existence of transverse heteroclinic(homoclinic) orbits in the Poincar\'e map for a perturbed one-degree of freedom
Hamiltonian system by measuring the ``distance'' between the stable and
unstable manifolds associated with the saddle points. This method was
generalized by Holmes and Marsden to the case of perturbed two-degree of
freedom Hamiltonian systems when the phase space is a product of the dual of
a Lie algebra and a set of action-angle variables. We briefly recall below this result; see Holmes and Marsden \cite{H-M} for proofs.

The setting is the following. The phase space is the product of the dual
of a Lie algebra $\mathfrak{g}$ and $\mathbb{R}^{2}$. The Hamiltonian has the form%
\begin{align}
H^{\varepsilon }(\mu ,\theta ,I)&
=F(\mu )+G(I)+\varepsilon H^{1}(\mu ,\theta,I)
+O(\varepsilon ^{2}) \nonumber \\
&= H^0(\mu, I) +  \varepsilon H^{1}(\mu ,\theta,I)
+O(\varepsilon ^{2}) 
\label{sisM}
\end{align}
where $\mu =(\mu _{1},...,\mu _{m})\in \mathfrak{g}^{\ast }$ and ($\theta ,I$) are coordinates on $\mathbb{R}^{2}$, with $\theta $ a $2\pi$-periodic
variable. It is also assumed that the Lie-Poisson system whose Hamiltonian is $F$ has a
heteroclinic (or homoclinic) orbit 
$\overset{\backsim }{\mu }(t)\in 
\mathfrak{g}^{\ast }.$ The oscillator frequency
\begin{equation*}
\Omega (I):=\frac{\partial G}{\partial I}
\end{equation*}
is assumed to be positive. The result is the following:

\begin{theorem}[Holmes-Marsden]
Suppose $\overset{\backsim }{\mu }(t)$ is a heteroclinic (or homoclinic)
orbit for the Lie-Poisson system whose Hamiltonian is $F,$ which lies in a two dimensional coadjoint orbit in $\mathfrak{g}^{\ast }$. Let $\overset{\backsim }{h}=F(\overset{\backsim }{\mu })$ be
the energy of the heteroclinic orbit and let $h>\overset{\backsim }{h}$ and 
$l^{0}=G^{-1}(h-\overset{\backsim }{h})$ be constants. Let $\{F,H^{1}\}(t,\theta ^{0})$ denote the Lie-Poisson bracket of $F(\mu )$ and $H^{1}(\mu ,\Omega (l^{0})t
+\theta ^{0},l^{0})$ evaluated at 
$\overset{\backsim }{\mu }(t)$. Let 
\begin{equation*}
M(\theta ^{0})=\frac{1}{\Omega (l^{0})}\overset{\infty }{\underset{-\infty }{\int }}\{F,H^{1}\}(t,\theta ^{0})dt
\end{equation*}
and assume $M(\theta ^{0})$ has simple zeros. Then for $\varepsilon >0$
sufficiently small, the Hamiltonian system (\ref{sisM}) contains transverse
heteroclinic orbits and hence Smale horseshoes on the energy surface $H^{\varepsilon }=h.$
\end{theorem}

Now we will prove that the system (\ref{sis2}) verifies the conditions of
the above theorem. The unperturbed system of (\ref{sis2}) on $\mathfrak{se}^{\ast }(2)$
has unstable critical points $(0,\pm M,0)$ lying on the $2$-dimensional
coadjoint orbit given by the cylinder
\begin{equation*}
\{(\mu _{1},\mu _{2},\mu _{3})\in \mathfrak{se}^{\ast }(2):\mu _{1}^{2}+\mu
_{2}^{2}=M^{2}\}.
\end{equation*}
On this coadjoint orbit we have the heteroclinic orbits given by
\begin{equation*}
\left\{ 
\begin{array}{c}
\mu _{1}(t)=\pm M\operatorname{sech}(Mt) \\ 
\mu _{2}(t)=\pm M\tanh (Mt) \\ 
\mu _{3}(t)=\pm M\operatorname{sech} (Mt)
\end{array}
\right.
\end{equation*}
that link the unstable critical points $(0,M,0)$ and $(0,-M,0)$.

On $\mathbb{R}^{2}$ the unperturbed system of (\ref{sis2}) is completely
integrable and in action-angle coordinates $(I, \theta)$, $u_1 = \sqrt{2I} \cos \theta$, $u_2 = \sqrt{2I} \sin \theta$, takes the form
\begin{eqnarray*}
&&%
\begin{array}{cc}
\overset{\cdot }{I}= & 0
\end{array}
\\
&&%
\begin{array}{cc}
\overset{\cdot }{\theta }= & 1.
\end{array}%
\end{eqnarray*}

Now writing the Hamiltonian (\ref{hap}) in the form
\begin{equation*}
H^{\varepsilon }(\mu _{1},\mu _{2},\mu _{3},I,\theta )
=\frac{1}{2}\left( \mu_{1}^{2}+2\mu _{2}^{2}
+\mu _{3}^{2}\right) 
+I -\varepsilon \mu_{3}\sqrt{2I}\sin \theta 
+O(\varepsilon ^{2}).
\end{equation*}
we are in the setting of the Holmes-Marsden theorem and we can write the
heteroclinic orbits for an energy level
\begin{equation*}
H^{0}=h=M^{2}+k,
\end{equation*}
where $M^{2}=\mu _{1}^{2}+\mu _{2}^{2}$ and $k$ is a constant, as
\begin{equation*}
\left\{ 
\begin{array}{l}
\mu _{1}(t)=\pm M\operatorname{sech}(Mt) \\ 
\mu _{2}(t)=\pm M\tanh (Mt) \\ 
\mu _{3}(t)=\pm M\operatorname{sech}(Mt) \\ 
\quad \; \; \, I=k \\ 
\quad \; \; \, \theta =t+\theta ^{0}.
\end{array}
\right.
\end{equation*}

The Melnikov function is:
\begin{eqnarray*}
M(\theta ^{0}) &=&-\int_{-\infty }^{\infty }
\sqrt{2k}M^{2}\operatorname{sech}(Mt)\tanh (Mt)
\sin (t+\theta ^{0})dt \\
&=&-\sqrt{2k}M^{2}\left( \int_{-\infty }^{\infty }\sinh (Mt)\cosh^{-2}(Mt)\sin (t)dt\right)\cos \theta ^{0} \\
&=&- \pi \sqrt{2k}  \operatorname{sech}\left(\frac{\pi }{2M}\right)\cos \theta ^{0},
\end{eqnarray*}
which has simple zeros as a function of $\theta^0$ and therefore the Hamiltonian system (\ref{sis2}) for 
$\varepsilon >0$ sufficiently small, has transverse heteroclinic orbits and
hence Smale horseshoes in a suitably chosen cross section of the constant
energy surface with $k>0.$

Acknowledgments. This research was partially
supported by the European Commission and the Swiss Federal
Government through funding for the Research Training Network
\emph{Mechanics and Symmetry in Europe} (MASIE) as well as
the Swiss National Science Foundation.


\begin{thebibliography}{9}
\bibitem{bok} O. Bokhove, On Hamiltonian balanced models. Preprint, Ninth 
\emph{Conf. on Atmospheric and Oceanic Waves and Stability}, San Antonio,
TX, Amer. Meteor. Soc., 367-368 (1993)

\bibitem{os} O. Bokhove and T. G. Shepherd, On Hamiltonian Balanced Dynamics
and the Slowest Invariant Manifold. \emph{J. Atmos. Sci., }\textbf{53, }%
276-297 (1996)

\bibitem{cam} R. Camassa, On the geometry of an atmospheric slow manifold. 
\emph{Physica D}, \textbf{84}, 357-397 (1995)

\bibitem{H-M} P.J. Holmes and J.E. Marsden, Horseshoe and Arnold diffusion
for Hamiltonian systems on Lie groups. \emph{Indiana University Math. J.} 
\textbf{32}, 273-309 (1983)

\bibitem{lor} E. N. Lorenz, On the existence of a slow manifold. \emph{J.
Atmos. Sci.},\textbf{\ 43}, 1547-1557 (1986)

\bibitem{mel} V. K. Melnikov, On the stability of the center for time
periodic perturbations. \emph{Trans. Moscow Math.}, \textbf{12}, 1-57 (1963)

\bibitem{wig} S. Wiggins, Global Bifurcations and Chaos, 
Springer-Verlag, New York, NY, 1988.
\end{thebibliography}
\end{document}